\documentclass[12pt]{article}
\newcommand{\vs}{\vspace{.1in}}
\newcommand{\ds}{\displaystyle}
\begin{document}
\title{The numbers behind Plimpton 322\footnote{The original
version of this article was dated August 6, 2011. It is
reworked here to 
take advantage of the new research
on the topic published by Britton, Proust and Shnider [2011]
later that year.}}
\author{Anthony Phillips\\Mathematics Department, Stony Brook University}

\date{\today}
\maketitle
\begin{abstract}

A mathematically and culturally natural 
modification of Evart Bruins' explanation of
the genesis of the numbers on the Old Babylonian
tablet Plimpton 322 gives
an economical accounting for the ``missing
pairs'' in his reconstruction. 
When the new
scheme is used to predict the numbers that would
follow those on Plimpton 322, 
the results add five new rows to those
listed by Price and Friberg in their hypothetical extension
of the content of Plimpton 322 to rows covering the edges and back.
\end{abstract}

\section{Introduction}
Plimpton 322 is probably the best-known Babylonian mathematical text. It
must be one of the most controversial; 
substantially divergent explanations have
been offered by very reputable scholars both for its purpose (why?)
and for its generation (how?). 
A recent article by Britton, Proust and Shnider [2011] includes
a comprehensive review and analysis of the interpretative efforts
to date, and can be taken as a ``state of the art'' document,
from the assyriological point of view. But Plimpton 322 is
(as these authors substantiate) a mathematical document and
as such has been subjected to mathematical analysis since its
first publication. The purpose of this note is to review
the purely mathematical aspects of this analysis, and to
suggest a minor and, one can argue, culturally natural
modification to the current understanding of the process by which the
particular numbers, which appear in the 15 rows of the tablet, were chosen.
When this new
criterion is used to extrapolate
beyond those 15 rows, it adds five new rows to those
that have been previously suggested, by Price [1964] and by Friberg [1981].

\section{The numbers on the tablet}

Here is a transcription\footnote{A recent high-quality photograph
of the tablet, front and back, is available on the website of the
Institute for the Study of the Ancient World:
http://isaw.nyu.edu/exhibitions/before-pythagoras/items/plimpton-322/}
   of the numerical part of the tablet, 
reproduced from [Robson 2001]. 
\begin{verbatim}
[1 59] 00 15                 1 59        2 49     KI.1 
[1 56 56] 58 14 50 06 15    56 07        1 20 25  KI.2 
[1 55 07] 41 15 33 45        1 16 41     1 50 49  KI.3 
1 53 10 29 32 52 16          3 31 49     5 09 01  KI.4 
1 48 54 01 40                1 05        1 37     KI.[5] 
1 47 06 41 40                5 19        8 01     [KI.6] 
1 43 11 56 28 26 40         38 11       59 01     KI.7 
1 41 33 45 14 3 45          13 19       20 49     KI.8 
1 38 33 36 36                8 01       12 49     KI.9 
1 35 10 02 28 27 24 26 40    1 22 41     2 16 01  KI.10 
1 33 45                     45           1 15     KI.11 
1 29 21 54 2 15             27 59       48 49     KI.12 
1 27 00 03 45                2 41        4 49     KI.13 
1 25 48 51 35 6 40          29 31       53 49     KI.14 
1 23 13 46 40               28          53        KI.15
\end{verbatim}\vs

Text in [brackets] is a reconstitution of missing material. 
The initial 1s are now generally accepted [Britton {\it et al.} 2011];
they occur along the cleavage between Plimpton 322 and a missing
left-hand portion of the tablet. 
This transcription includes six generally accepted
corrections of the cuneiform entries; again, see [Britton {\it et al.} 2011],
although Friberg [2007] maintains 3 13 in Row 2, Column 3 and 56
in row 3, Column 2.
\vs

The four columns of the tablet
 may be labeled from left to right $A$, $S$, $D$ and $N$, the last
a row number running from 1 to 15. 
\vs

The uncontested (except for Friberg [2007]) mathematical properties of the 
numbers on the tablet are as follows 
\begin{enumerate}
\item  The numbers in column $A$ decrease, fairly regularly, from top to bottom.
\item  Each number in column $A$ is a perfect square and one more than a perfect square.

\item  In each row, $D^2-S^2$ is a perfect square.
\item In each row, $A = D^2/(D^2-S^2)$. For example, the fourth row contains the sexagesimal numbers
$$\begin{array}{cccc}
A                    &   S~~  &       D~~~  &       N\\
1~ 53~ 10~ 29~ 32~ 52~ 16~~~~ &   3 ~31 ~49~~~~ &  5 ~09 ~01~~~~&   4 \\
\end{array}$$
(here $D^2-S^2 = 3~ 45~00^2$)
and $D^2/(D^2-S^2)
 = 1~ 53~ 10 ~29 ~32 ~52 ~16 = A$. 
\item  In each row, the numbers in column $S$ and column $D$ have no common factors. (One exception: row 11).
\end{enumerate}
Additionally, the column headings for columns $S$ and $D$ are read to refer to the ``short-side'' or ``front''
and the ``diagonal,'' presumably of a rectangle or a right triangle.

\section{Where do these numbers come from?}

The standard explanation dates back to the publication of
the tablet by Neugebauer and Sachs [1945]
and has persisted, with some elaboration, until today ([Britton {\it et al}, 
2011]). 

Their explanation involves a virtual column $L$ 
(not on the table, although they speculate that it
might have been on the missing half), related to
the $S$ and $D$ columns by $L^2=D^2-S^2$, so that
$L,S,D$ (long, short, diagonal) form a Pythagorean
triple in modern terminology; this also means that on the
tablet $A = (D/L)^2$. 

In those terms, here 
is how Neugebauer and Sachs understood the problem: ``How were the
mathematicians of the Old Babylonian period able not only to solve
the Pythagorean equation ... in integers but to adapt the solutions
to the further condition that the proportion $D/L$ decrease from step
to step by a number deviating very little from one-sixtieth?''

It has been known at least since Euclid (X, Proposition 29, Lemma 1)
that (in the reformulation used by Neugebauer and Sachs)
if $P$ and $Q$ are numbers satisfying ($i$) $(P,Q)=1$,
($ii$) $P$ and $Q$ not both odd, ($iii$) $P>Q$,
then the numbers 
$$2PQ, P^2-Q^2, P^2+Q^2~~~~~(*)$$ form a primitive
Pythagorean triple, i.e. one with no shared factor among
the three. Furthermore, since in any primitive Pythagorean
triple $a^2+b^2=c^2$ the numbers $a$ and $b$ must have
opposite parity, one can relabel the odd one as $o$ and
the even\footnote{In fact all the entries
in the $S$ column on Plimpton 322 are odd except for row 15; 
so $P=\sqrt{(D+S)/2}$, etc., except that in (the corrected)
row 15 $P=\sqrt{53+45}/2$.}
one as $e$; 
then $P$ and $Q$ can be retrieved from
$a, b, c$ by $P=\sqrt{(c+o)/2}$, $Q=\sqrt{(c-o)/2}$. One
can check that $(P,Q)=1$ and that the two cannot both be odd.
So there is a a one-one correspondence
between primitive Pythagorean triples $a, b, c$ and pairs $P,Q$
satisfying ($i$), ($ii$), ($iii$). 

Neugebauer and Sachs assume that the Old Babylonians were aware
of this fact, and give
the $P, Q$ pairs corresponding to the columns of Plimpton 322:
$$\begin{array}{crrc|crrc|crr}
\mbox{row}&P&Q&~~&\mbox{row}&P&Q&~~&\mbox{row}&P&Q\\
1&12 &5 &~~&2 &1,4 &27 &~~&3 &1,15 &32 \\
4&2,5 &54 &~~&5 &9 &4 &~~&6 &20 &9 \\
7&54 &25 &~~&8 &32 &15 &~~&9 &25 &12 \\
10&1,21 &40 &~~&11 &2 &1 &~~&12 &48 &25 \\
13&15 &8 &~~&14 &50 &27 &~~&15 &9 &5
\end{array}$$

There are a couple of minor problems, since
in row 11 $S$ and $D$ are not relatively prime
(the theorem would give $S=3, D=5$) and in row
15 they correct the tablet's $S=56, D=53$ (obviously
an error since that would make the short side longer than
the diagonal) to $S=56, D=1, 46$, not
relatively prime, leading to $P=9, Q=5$, both odd. 
The correction generally accepted now is $S=28, D=53$ which
gives $P=7, Q=2$.

These problems do not seem significant; but
the question remains, where did these $P,Q$ pairs
come from. Neugebauer and Sachs first
remark: ``With the single exception of $2,5 = 5^3$
all the numbers $P$ and $Q$ of our list are contained in
the group of regular numbers which constitute the
`reciprocal tables' ''  

They are referring to the standard reciprocal table,
for which [Neugebauer and Sachs, 1945] lists 13 specimen tablets. The table gives
some 30 ``standard'' pairs of reciprocal numbers: numbers
which multiply to a power of 60, i.e. to 1 in the
sexagesimal floating-point notation of the period.\vs

The currently accepted (but not in [Friberg, 2007])
reading for row 15, which leads to
$P=7, Q=2$, is incompatible with this analysis: their exceptional
$2,5$ is still a regular number but $7$ is definitely not.\vs

On the other hand, Neugebauer and Sachs point the way to a
different approach to the problem by remarking that 
$$A = (D/L)^2 = [\frac{\ds P^2+Q^2}{\ds 2PQ}]^2
=[\frac{1}{2}(P\overline{Q}+Q\overline{P})]^2),$$
including
$$D/L = \frac{1}{2}(P\overline{Q}+Q\overline{P})~~~~~~(**),$$ 
where $\overline{X}$
represents the reciprocal of $X$, and highlighting the
reciprocal pairs $P\overline{Q}, Q\overline{P}$ in their
summary:

``Our final result, then, is that our tablet was calculated by
selecting numbers $P\overline{Q}$ and $Q\overline{P}$ from combined
multiplication tables such that $(**)$ has a value as near as possible
to the required [for uniform linear decrease] values of $D/L$;
Pythagorean numbers were then formed with these values of $P$
and $Q$ according to $(*)$.''\vs

Neugebauer and Sachs go so far as to suggest that an alternative formulation might
be considered, ``using one parameter $\alpha$ and its reciprocal 
$\overline{\alpha}$, where $\alpha=\frac{P}{Q}$'' and they show 
what this would give for the
first four lines of the tablet:
$$\begin{array}{l|r|l|l}
P & Q & \alpha = P\overline{Q}&\overline{\alpha}=Q\overline{P}\\\hline
12&5&2;24&0;25\\
1,4&27&2;22,13,20&0;25,18,45\\
1,15&32&2;20,37,30&0;25,36\\
2,5&54&2;18,53,20&0;25,55,12
\end{array}$$

The computation would presumably continue:
$$\begin{array}{l|l|l|l|l|l}
\alpha = P\overline{Q}&\overline{\alpha}=Q\overline{P}&S\overline{L}=\frac{1}{2}( \alpha- \overline{\alpha})& D\overline{L}=\frac{1}{2}( \alpha+ \overline{\alpha})&S&D\\\hline 
2;24      &0;25      &0;59,30       &1;24,30        & &\\  
2;22,13,20&0;25,18,45&0;58,27,17,30 &1;23,46,2,30  & &\\ 
2;20,37,30&0;25,36&0;57,30,45&1;23,6,45  & &\\
2;18,53,20&0;25,55,12 &0;56,29,4 &1;22,24,16 & & 
\end{array}$$
where $S$ and $D$, being regular and relatively prime, can be
obtained from $S\overline{L}$ and $D\overline{L}$ by casting
out common factors of 2, 3, and 5.\vs

But they dismiss this formulation in favor of the ``simple numbers
$P$ and $Q$'' and equation $(*)$.\vs

Price [1964] put the selection of $P$ and $Q$ foremost. He
remarked that the values found by Neugebauer and Sachs satisfy
$$ 1<Q<60 ~~~~~~f<P/Q<g ~~~~~~(***)$$
``where $f$ is such that $9/5=1,48$ is included but $16/9=1,46,10$
is not, and $g$ is such that $12/5=2,25$ is included but $5/2=2,30$
is not''  and that ``the tablet contains [entries hypothetically
corresponding to] all and only all those pairs of regular numbers,
prime to each other, and satisfying the inequalities $(***)$.''\vs

The words ``all and only all'' are important because they must
be part of a coherent mathematical explanation. Now the
apparent arbitrariness in the selection of $P$ and $Q$ has
been considerably reduced, to the choice of the bounds on $Q$ and on $P/Q$.
The $1<Q<60$ can be accepted as ``natural'' since it
defines single (sexagesimal) place integers larger than 1. 
The bounds $f$ and $g$ are less easy. For $g$ Price continues
with the interpretation of Plimpton 322 as a list of Pythagorean
triples: ``it has already been
remarked by Neugebauer that the upper bound, $g$, is certainly intended
to correspond with the isosceles right-angle triangle. Indeed, the
first parameter on the list, (12:5) gives the triangle 2,0/1,59/2,49
where the sides come closest to equality.'' For $f$ Price states that
``it is easy to see that a different limit of $f=1$ would be very natural.''
Staying with Pythagorean triangles, it is the limit as one of the
acute angles goes to 0. This new limit would allow an additional
23 $(P,Q)$ pairs; Price joins this fact with the observation that the column
dividers on the tablet are repeated on the back (so there would
be ``almost exactly the right amount of space for the complete array
to be given in order''\footnote{Note that about 
5 extra rows would need to fit on the top and the bottom
of the tablet; but the column divisions were not carried over to
those surfaces ([Friberg 1981], Fig.\ 1.3).}.)
Price's summation: ``It appears then that the Plimpton tablet
is based upon a complete collection of all the Pythagorean triples that can
be produced from pairs of parameters that are both regular and prime to
each other. The first of the parameters is chosen so that it is a single-place
sexagesimal integer (other than unity) and the second is chosen so that
the sides of the Pythagorean triangle satisfy the inequality
$0<P^2-Q^2<2PQ$.'' This last inequality is equivalent to 
$1 < P/Q < 1+\sqrt{2}$.\vs

Friberg [1981] took up Price's analysis and extended the $(P,Q)$
collection to include ``all pairs $(P,Q)$ with $0<P<1~ 00$,
$Q\leq 2~ 05$ and $0.15 \leq Q/P < 1$,'' corresponding to
$1<P/Q<5.6$, and generating 22 additional
rows which he numbered $-1,\dots, -22$. In his Figure 2.2 he
lists all 60 pairs, indexed by their corresponding $(P,Q)$;
 but he keeps the alternative ``reciprocal pairs'' explanation
alive by giving $P\overline{Q}$
and $Q\overline{P}$ for each item. His proposed reconstruction of
Plimpton 322 (his Figure 2.3) adds two new columns:
$\frac{1}{2}(P\overline{Q}-Q\overline{P})$ and
$\frac{1}{2}(P\overline{Q}+Q\overline{P})$. When he returns
to the question in [2007] his reconstruction has yet two more columns
on the left: $P\overline{Q}$ and $Q\overline{P}$ (Fig. A8.3).     \vs

Britton, Proust and Shnider [2011] continue this tradition:
``In a sense, the `point of departure' for the construction of
the text was indeed the simple regular numbers $P$ and $Q$,
as suggested by Neugebauer and Sachs, but with $Q$ limited to
sexagesimal digits, and $P$ constrained only by the condition that
$P/Q$ be greater than 1 and less than $\sqrt{2}+1$'' 
(they do not consider Friberg's examples with $\sqrt{2}+1
<P/Q<5.6$, corresponding to his negatively numbered rows,
presumably on the grounds that the ``short side''
$S$ thus defined would be longer than the other side,
$\sqrt{D^2-S^2}$). But they
argue convincingly, through a detailed analysis of the four
non-transcriptional errors in the tablet, that $P$ and $Q$
enter into the calculation of the entries in the tablet only
through their repackaging as $\frac{1}{2}(P\overline{Q}-Q\overline{P})$ and
$\frac{1}{2}(P\overline{Q}+Q\overline{P})$.
\vs

Their analysis, as well as Friberg's choice [2007] for two of the 
missing columns, suggests reinvestigating Neugebauer and Sachs'
rejected  hypothesis, that a list of reciprocal
pairs $\alpha, \overline{\alpha}$ underlies the set of numbers on
Plimpton 322, and that $P$ and $Q$ are artifacts of the analysis,
satisfying $P\overline{Q}=\alpha$. 
This alternative route was proposed by Bruins [1949], by Buck [1980]
citing an unpublished article by D.\ L.\ Voils, by Schmidt [1980]
and most recently by Robson [2001]. The only real objection to
this solution is that it requires the author of
Plimpton 322 to have had access to a table containing
4-place sexagesimal reciprocals. Robson [2001] addresses this problem
at length;
whether that analysis is satisfactory is cannot appropriately
be discussed here, but the purpose of this article is to present
a mathematical detail supporting the Bruins-Robson approach.

\section{A missing element in the reciprocal-pairs explanation}  

Let us suppose then that the 15 rows on Plimpton 322 were generated
from 15 elements of a list of $\alpha, \overline{\alpha}$ reciprocal pairs. 
Then two questions arise: what kind of list was it, and what determined
the selection of those particular 15?\vs

The 15 $\alpha, \overline{\alpha}$ pairs, starting with the four
given by Neugebauer and Sachs, are ([Bruins, 1957], 
Table II; [Robson, 2001], Table 9; [Friberg, 2007],
Fig. A8.3):
$$\begin{array}{cll|cll|cll}

\mbox{Row} &    \alpha   &           \overline{\alpha}&
\mbox{Row} &    \alpha   &           \overline{\alpha}&
\mbox{Row} &    \alpha   &           \overline{\alpha}\\
1 &   2~ 24         &    25~	  &
2 &   2~ 22~ 13~ 20 &    25~ 18~ 45&
3 &   2~ 20~ 37~ 30 &    25~ 36	\\
4 &   2~ 18~ 53~ 20 &    25~ 55~ 12&
5 &   2~ 15         &    26~ 40	&
6 &   2~ 13~ 20     &    27	\\
7 &   2~ 09~ 36  &   27~ 46~ 40	&
8 &   2~ 08         &    28~ 07~ 30&
9 &   2~ 05         &    28~ 48	\\
10&   2~ 01~ 30   &    29~ 37~ 46~ 40&
11&   2~ 00         &    30	&
12&   1~ 55~ 12      &    31~ 15\\
13&   1~ 52~ 30     &    32	&
14&   1~ 51~ 06~ 40 &    32~ 24	&
15&   1~ 48         &    33~ 20	
\end{array}$$
These generate the $S$ and $D$ columns as initially suggested
by Neugebauer and Sachs, {\it via} $\frac{1}{2}(\alpha - \overline{\alpha})$,
$\frac{1}{2}(\alpha + \overline{\alpha})$ and removing common factors,
as in the examples above,
with the usual exception for row 11 where the common factor 15
is left in.\vs

Three remarks are in order.
\begin{enumerate}
\item The $\alpha$ decrease monotonically from row 1 to 15 (as
opposed to the chaotic behavior of the $P,Q$ pairs).
\item Every $\alpha$ and $\overline{\alpha}$ is this list has
at most four sexagesimal places, and the fourth place is 0 or
a multiple of 10.
\item There are exactly fifteen such pairs of ``round-number"
reciprocals in the range between $2~ 24 $ and $1~ 48$, inclusive.
\end{enumerate}

I do not believe that last remark has been made before. But the
``multiple of 10'' criterion also
gives an elementary explanation for why the six four-place
pairs listed by Bruins [1949, Table II] and Robson [2001, Table 8], 
while in the same 
range, do not contribute to the table. \vs

These criteria would define an ordered list of reciprocal pairs, from
which 15 consecutive elements can be used to generate the numbers
on Plimpton 322. Britton {\it et al.}'s [2011] geometric bounds 
$1<\alpha<\sqrt{2}+1$, necessary for congruence with the
column headings ``short side'' and ``diagonal,'' now imply a
different continuation of the series beyond the tablet's row 15.

\section{Predictions for additional rows}
 If an ordered list 
of 4-place reciprocal pairs $\alpha,\bar{\alpha}$
divisible by 10, with $1<\alpha<\sqrt{2}+1$, was in fact the primary 
given in constructing Plimpton 322, then that list should suggest
how the tablet was meant to be extended. Comparing those prodictions 
with the
hypothetical extensions of the Plimpton list 
 in [Price 1964] and [Friberg 1981] shows
substantial, but not perfect, agreement. 

Here
is a comparison of the reconstructions; the row numbers
given are from the Friberg and Price/Friberg lists; new
elements are identified as {\it i}, {\it ii}, etc.

$$\begin{array}{clrrrrlrrrr}
\mbox{No.}& ~    &~ & \alpha &~  &~  & ~     &~  & \overline{\alpha}&~  &~  \\
   16&~~~~~~~&    1& 46& 40&  0&~~~~~~~&        33& 45&  0&  0\\
   i&~~~~~~~&    1& 44& 10&  0&~~~~~~~&        34& 33& 36&  0\\
   ii&~~~~~~~&    1& 42& 24&  0&~~~~~~~&        35& 09& 22& 30\\
   17&~~~~~~~&    1& 41& 15&  0&~~~~~~~&        35& 33& 20&  0\\
   18&~~~~~~~&    1& 40&  0&  0&~~~~~~~&        36&  0&  0&  0\\
   19&~~~~~~~&    1& 37& 12&  0&~~~~~~~&        37& 02& 13& 20\\
   20&~~~~~~~&    1& 36&  0&  0&~~~~~~~&        37& 30&  0&  0\\
   21&~~~~~~~&    1& 33& 45&  0&~~~~~~~&        38& 24&  0&  0\\
   22&~~~~~~~&    1& 30&  0&  0&~~~~~~~&        40&  0&  0&  0\\
   23&~~~~~~~&    1& 28& 53& 20&~~~~~~~&        40& 30&  0&  0\\
   24&~~~~~~~&    1& 26& 24&  0&~~~~~~~&        41& 40&  0&  0\\
   25&~~~~~~~&    1& 25& 20&  0&~~~~~~~&        42& 11& 15&  0\\
   26&~~~~~~~&    1& 24& 22& 30&~~~~~~~&        42& 40&  0&  0\\
   27&~~~~~~~&    1& 23& 20&  0&~~~~~~~&        43& 12&  0&  0\\
   28&~~~~~~~&    1& 21&  0&  0&~~~~~~~&        44& 26& 40&  0\\
   29&~~~~~~~&    1& 20&  0&  0&~~~~~~~&        45&  0&  0&  0\\
   iii&~~~~~~~&    1& 18& 07& 30&~~~~~~~&        46& 04& 48&  0\\
   30&~~~~~~~&    1& 16& 48&  0&~~~~~~~&        46& 52& 30&  0\\
   31&~~~~~~~&    1& 15&  0&  0&~~~~~~~&        48&  0&  0&  0\\
   32&~~~~~~~&    1& 12&  0&  0&~~~~~~~&        50&  0&  0&  0\\
   33&~~~~~~~& 1&11&64&0 &~~~~~~~& 50&37&30&0\\
  iv&~~~~~~~& 1&09&26&40 &~~~~~~~& 51&50&24&0\\
   34&~~~~~~~& 1&07&30&0 &~~~~~~~&  53&20&0&0\\
   35&~~~~~~~& 1&06&40&0 &~~~~~~~&  54&0&0&0\\
   36&~~~~~~~& 1&04&48&0 &~~~~~~~&  55&33&20&0\\
   37&~~~~~~~& 1&04&0&0  &~~~~~~~& 56&15&0&0\\
   38&~~~~~~~& 1&02&30&0 &~~~~~~~&  57&36&0&0\\
 v&~~~~~~~& 1&00&45&0 &~~~~~~~&  59&15&33&20\\
\end{array}$$

The new pairs $i, \dots, v$, seem just as natural as the
ones generating the  Price/Friberg candidates for the ``left out'' rows.
$i, iii, iv$ and $v$ are related by  $(1/2,2)$ to Friberg's \# $-17$,  
\# $-4$ and to rows 4, 10 on Plimpton 322, respectively, while
$ii$ is $(9,1/9)$ times the dual of the pair in Friberg's \# $-21$.
They all appear in the Seleucid tablet AO 6456
[Neugebauer 1935]. On the other hand adding five extra rows to the
list makes it very unlikely that the top, bottom and back of the
tablet could have been planned to contain them all. \vs

\noindent
{\Large \bf References}\vs

\noindent
Britton, J., C. Proust and S. Shnider, Plimpton 322, a review and
a different prespective, {\it Archive for History of Exact Sciences}
{\bf 65} (2011) 519-566.\vs

\noindent
Bruins, E.\ M., On Plimpton 322, Pythagorean numbers in Babylonian
Mathematics, {\it K.\ Ned.\ Akad.\ van Wet.
Proceedings} {\bf 52} (1949) 629-632.\vs


\noindent
Bruins, E.\ M., Pythagorean triads in Babylonian Mathematics,
{\it Math.\ Gaz.}, {\bf 41} (1957) 25-28.\vs

\noindent
Buck, R. Creighton, Sherlock Holmes in Babylon, {\it Am.\ Math.\ Mon.}
{\bf 57} (1980) 335-345.\vs
 
\noindent 
Friberg, J\"{o}ran, Methods and Traditions of Babylonian Mathematics:
Plimpton 322, Pythagorean triples, and the Babylonian Triangle
Parameter Equations, {\it Hist,\ Math.} {\bf 8} (1981)
277-318. \vs

\noindent
Friberg, J\"{o}ran, {\it A Remarkable Collection of Babylonian Mathematical 
Texts,} Springer, New York, 2007.\vs

\noindent
Neugebauer, Otto, {\it Mathematische Keilschrifttexte, I,} Springer,
Berlin, 1935.\vs 

\noindent
Neugebauer, Otto and A. Sachs, {\it Mathematical Cuneiform Texts,} American
Oriental Society, New Haven, 1945.\vs


\noindent 
Price, Derek J.\ de Solla, The Babylonian ``Pythagorean Triangle" Tablet,
{\it Centaurus} {\bf 10} (1964) 219-231.\vs

\noindent
Robson, Eleanor, Neither Sherlock Holmes nor Babylon: A Reassessment
of Plimpton 322, {\it Hist.\ Math.} {\bf 28} (2001) 167-206. \vs


\noindent
Schmidt, Olaf, On Plimpton 322. Pythagorean Numbers in Babylonian
Mathematics, {\it Centaurus} {\bf 24} (1980) 4-13.
\end{document}